\documentclass[14pt,a4paper]{amsart}
\usepackage[14pt]{extsizes}
\usepackage{relsize}
\usepackage{amsmath}
\usepackage{amssymb}
\usepackage{tikz-cd}
\usepackage{tikz}
\usepackage{amscd}
\usepackage[left=3cm,right=2cm,top=3cm,bottom=2cm]{geometry}
\usepackage{abstract}

\newtheorem{Theorem}{Theorem}
\newtheorem{Lemma}{Lemma}

\newtheorem{Definition}{Definition}
\newtheorem{Remark}{Remark}

\title{New series of moduli components of rank 2 semistable sheaves on $\mathbb{P}^{3}$ with singularities of mixed dimension}

\author[A. Ivanov]{Aleksei Ivanov}
\address{Department of Mathematics\\
National Research University
Higher School of Economics\\
6 Usacheva Street\\ 
119048 Moscow, Russia}
\email{anivanov$_{-}$1@edu.hse.ru}

\begin{document}

\maketitle

\begin{abstract}
\noindent We construct a new infinite series of irreducible components of the Gieseker-Maruyama moduli scheme $\mathcal{M}(k), ~ k \geq 3$ of coherent semistable rank 2 sheaves with Chern classes $c_1=0,~ c_2=k,~ c_3=0$ on $\mathbb{P}^3$ whose general points are sheaves with singularities of mixed dimension. These sheaves are constructed by elementary transformations of stable and properly $\mu$-semistable reflexive sheaves along disjoint union of collections of points and smooth irreducible curves which are rational or complete intersection curves. As a special member of this series we obtain a new component of $\mathcal{M}(3)$.


\end{abstract}

\section{Introduction}

Let $\mathcal{M}(0,k,2n)$ be the Gieseker-Maruyama moduli scheme of semistable rank-2 sheaves with Chern classes $c_1=0,\ c_2=k,\ c_3=2n$ on the projective space $\mathbb{P}^3$. Denote $\mathcal{M}(k)=\mathcal{M}(0,k,0)$. By the singular locus of a given $\mathcal{O}_{\mathbb{P}^3}$-sheaf $E$ we understand the set $\mathrm{Sing}(E)=\{x\in\mathbb{P}^3\ |\ E$ is not locally free at the point $x\}$. $\mathrm{Sing}(E)$ is always a proper closed subset of $\mathbb{P}^3$ and, moreover, if $E$ is a semistable sheaf of nonzero rank, every irreducible component of $\mathrm{Sing}(E)$ has dimension at most 1. For simplicity we will not make a distinction between a stable sheaf $E$ and corresponding isomorphism class $[E]$ as a point of moduli scheme. Also by a general point we understand a closed point belonging to some Zariski open dense subset.

Any semistable rank-2 sheaf $[E] \in \mathcal{M}(k)$ is torsion-free, so it fits into the exact triple
$$
0 \longrightarrow E \longrightarrow E^{\vee \vee} \longrightarrow Q \longrightarrow 0,
$$
where $E^{\vee \vee}$ is a reflexive hull of $E$ and $\text{dim}~Q \leq 1$. Conversely, take a reflexive sheaf $F$, a subscheme $X \subset \mathbb{P}^{3}$, an $\mathcal{O}_{X}$-sheaf $Q$ and a surjective morphism $\phi: F \twoheadrightarrow Q$, then one can show that the kernel sheaf $E:=\text{ker}~\phi$ is semistable when $F$ and $Q$ satisfy some mild conditions. We call a sheaf $E$ an \textit{elementary transform} of $F$ along $X$. In general, an elementary transform of a sheaf $F$ can be defined as follows.

\begin{Definition}
An elementary transform of a sheaf $F$ along an element $[F \overset{\phi}{\twoheadrightarrow} Q] \in \emph{Quot}^{P}(F)$ is a sheaf $E:=\emph{ker}~\phi$.
\end{Definition}

\noindent In fact, all known irreducible components of the moduli schemes $\mathcal{M}(k)$ general points of which correspond to non-locally free sheaves are constructed by using elementary transformations of stable reflexive sheaves.

More precisely, in \cite{JMT2} there were found two infinite series $\mathcal{T}(k,n)$ and $\mathcal{C}(d_1,d_2,k-d_1 d_2)$ of irreducible components of $\mathcal{M}(k)$ which (generically) parameterize stable sheaves with singularities of dimension 0 and pure dimension 1, respectively. General points of components of the first series are elementary transforms of stable reflexive sheaves along unions of $n$ distinct points in $\mathbb{P}^{3}$, while those of the second series are elementary transforms of instanton bundles of charge $k - d_1 d_2$ along smooth complete intersection curves of degree $d_1 d_2$.

Next, in \cite{IT} there were constructed three components of $\mathcal{M}(3)$ parameterizing sheaves with singularities of mixed dimension. General sheaves of these components are elementary transforms of stable reflexive sheaves with Chern classes $(c_2, c_3)=(2, 2), \ (2, 4)$ along a disjoint union of a projective line and a collection of points in $\mathbb{P}^{3}$. This approach was generalized in \cite{AJT} by doing elementary transformations of stable reflexive sheaves with other Chern classes along a disjoint union of a projective line and a collection of points in order to construct infinite series of components of $\mathcal{M}(-1,c_2,c_3)$.

Also it is worth to note that in \cite{JMT1} there were constructed the certain collections of divisors of the boundaries $\partial \mathcal{I}(k)=\overline{\mathcal{I}(k)} \setminus \mathcal{I}(k)$ of instanton components of $\mathcal{M}(k)$ for each $k$. General sheaves of these divisors are elementary transforms of instanton bundles along rational curves.

The present paper is devoted to further generalization of these results. Namely, we construct an infinite series of irreducible moduli components which includes the components parameterizing non-locally free sheaves constructed in \cite{JMT1, JMT2, IT} as special cases. Similar to the construction in \textit{loc. cit.}, the general sheaves $E$ of the new components are obtained by the elementary transformations of the following form
$$
0 \longrightarrow E \longrightarrow F \longrightarrow L \oplus \mathcal{O}_{W} \longrightarrow 0,
$$
where $F$ is a stable or properly $\mu$-semistable reflexive (non-locally free) rank-2 sheaf, $L$ is a line bundle on a smooth connected curve $C$ in $\mathbb{P}^{3}$ which is either rational or complete intersection curve, $W \subset \mathbb{P}^{3}$ is a collection of points. In order to simplify computations we require that $C \cap W = \emptyset$ and $\text{Sing}(F) \cap (C \sqcup W) = \emptyset$. One can show that the singularity set of the sheaf $E$ is $\text{Sing}(F) \sqcup C \sqcup W$, so it has mixed dimension. Moreover, $\text{Sing}(E)$ does not coincide with any other singularity set of the sheaves from the known components of $\mathcal{M}(k)$, so the components of the proposed series are really new.   

Since a complete enumeration of components of $\mathcal{M}(k)$ for small values of $k$ is of particular interest, it is worth to note that this series contains a new component of $\mathcal{M}(3)$. In short, the dense subset of this component can be obtained by doing elementary transformations of properly $\mu$-semistable reflexive rank-2 sheaves $F$ with $(c_1, c_2, c_3) = (0, 1, 2)$ along the sheaf $L = \mathcal{O}_{C}(2)$, where $C \subset \mathbb{P}^{3}$ is a smooth conic.

The paper is organized in the following way. In Section 2 we recall the necessary facts about moduli spaces of stable and $\mu$-semistable reflexive sheaves. Section 3 is devoted to the description of the new series of moduli components. Finally, in Section 4 we prove that the described components are irreducible.

\vspace{3mm}\noindent{\bf Acknowledgements.}
The work was supported in part by Young Russian Mathematics award and by the Russian Academic Excellence Project ‘‘5-100’’. I would like to thank A. S. Tikhomirov for usefull discussions and D. Markushevich for the opportunity to give a talk about results of this paper on the conference "Integrable systems and automorphic forms" (University of Lille-1, 2019). 

\section{Reflexive rank-2 sheaves}

The moduli scheme $\mathcal{R}(0,m,2n)$ parameterizing stable reflexive rank-2 sheaves on $\mathbb{P}^{3}$ with Chern classes $c_1=0, \ c_2=m, \ c_3=2n$ can be considered as open subset of the Gieseker-Maruyama moduli scheme $\mathcal{M}(0,m,2n)$, so it is a quasi-projective scheme (see \cite{SRS}). It is known that for $(m,n)=(2,1), \ (2,2), \ (3,4)$ this scheme is smooth, irreducible and rational; for $(m,n) = (3,2)$ it is irreducible and reduced at general point; for $(m,n)=(3, 1), \ (3,3)$ the corresponding reduced scheme is irreducible (see \cite{Chang}). 

In the paper \cite{JMT2} the infinite series of irreducible components $\mathcal{S}_{a,b,c}$ of the moduli schemes $\mathcal{R}(0,m,2n)$ is described. Sheaves from these components fit into the following exact triple
\begin{equation}\label{reflexive series}
0 \rightarrow a \cdot \mathcal{O}_{\mathbb{P}^3}(-3) \oplus b \cdot \mathcal{O}_{\mathbb{P}^3}(-2) \oplus c \cdot \mathcal{O}_{\mathbb{P}^3}(-1) \rightarrow (a+b+c+2) \cdot \mathcal{O}_{\mathbb{P}^3} \rightarrow F(k) \rightarrow 0,
\end{equation}
\noindent where $a,~ b,~ c$ are arbitrary non-negative integers such that $3a+2b+c$ is non-zero and even, $k:=\frac{3a+2b+c}{2}$.  The corresponding Chern classes of these sheaves can be expressed through the integers $a,~ b,~ c$ in the following way
\begin{equation}\label{pol m}
c_2(F)=\frac{1}{4}(3a+2b+c)^{2}+\frac{3}{2}(3a+2b+c)-(b+c),
\end{equation}
\begin{equation}\label{pol n}
c_3(F)=27 {a+2 \choose 3} + 8 {b+2 \choose 3} + {c+2 \choose 3} +
\end{equation}
$$
+ 3(3a + 2b +5)ab + \frac{3}{2}(2a+c+4)ac + (2b+3c+3)bc + 6abc.
$$
\noindent The components $\mathcal{S}_{a,b,c}$ are smooth. Moreover, they have expected dimension $8m-3$ which implies that $\text{Ext}^{2}(F,F) = 0$ for any sheaf $[F] \in \mathcal{S}_{a,b,c}$ (see \cite[Lemma 5]{JMT2}).

Also we can construct a scheme $\mathcal{V}(0,m,2n)$ parameterizing some reflexive properly $\mu$-semistable sheaves with the corresponding Chern classes in the following way. Consider the Hilbert scheme $\text{Hilb}_{m,g}(\mathbb{P}^{3})$ of smooth space curves of degree $m$ and genus $g$; let $n=g+2m-1$. Now denote by $\mathcal{Z} \hookrightarrow \text{Hilb}_{m,g}(\mathbb{P}^{3}) \times \mathbb{P}^{3}$ the corresponding universal curve and $\text{pr}: \text{Hilb}_{m,g}(\mathbb{P}^{3}) \times \mathbb{P}^{3} \longrightarrow \text{Hilb}_{m,g}(\mathbb{P}^{3})$ the projection onto the first factor. Also introduce the following definition 
\begin{Definition}
Let $S$ be a scheme. Let $\mathcal{E}$ be a coherent $\mathcal{O}_{S}$-sheaf. We denote by \emph{$\textbf{P}(\mathcal{E}):=\text{Proj}(\text{Sym}_{\mathcal{O}_{S}}(\mathcal{E}))$} the Proj construction of the sheaf of graded \emph{$\mathcal{O}_{S}$}-algebras \emph{$\text{Sym}_{\mathcal{O}_{S}}(\mathcal{E})$}.
\end{Definition}
\noindent We define the scheme $\mathcal{V}(0,m,2n)$ as an open subset of $\textbf{P}((\text{pr}_{*}\omega_{\mathcal{Z}}(4))^{\vee})$ the points $(Y, \ \mathbb{P}\xi) \in \textbf{P}((\text{pr}_{*}\omega_{\mathcal{Z}}(4))^{\vee})$ of which satisfy the following property
$$
\xi \in \text{H}^{0}(\omega_{Y}(4)) \ \text{generates} \ \omega_{Y}(4) \ \text{except at finitely many points}.
$$
By the construction we have the formula for the dimension of this scheme
\begin{equation}\label{dim_strictmu}
\text{dim}~\mathcal{V}(0,m,2n) = \text{dim}~\text{Hilb}_{m,g}(\mathbb{P}^{3}) + \text{dim}~\mathbb{P}(\text{H}^{0}(\omega_{Y}(4))) =
\end{equation}
$$
= h^{0}(N_{Y/\mathbb{P}^{3}}) + h^{0}(\omega_{Y}(4)) - 1,
$$
\noindent where $Y$ is an arbitrary curve from $\text{Hilb}_{m,g}(\mathbb{P}^{3})$. Next, note that due to the isomorphism $\text{H}^{0}(\omega_{Y}(4)) \simeq \text{Ext}^{1}(I_{Y}, \mathcal{O}_{\mathbb{P}^{3}})$ any point $(Y, \mathbb{P}\xi) \in \mathcal{V}(0,m,2n)$ uniquely defines the sheaf $F$ which satisfies the following exact triple
\begin{equation}\label{serre for non-stable}
0 \longrightarrow \mathcal{O}_{\mathbb{P}^{3}} \longrightarrow F \longrightarrow I_{Y} \longrightarrow 0.
\end{equation}
One can show that $F$ is a reflexive properly $\mu$-semistable rank-2 sheaf with Chern classes $c_1=0, \ c_2=m, \ c_3=2n$. Conversely, any such sheaf $F$ satisfies the triple above, so it determines the point of $\mathcal{V}(0,m,2n)$. Therefore, there exists one-to-one correspondence between points of $\mathcal{V}(0,m,2n)$ and some family of reflexive properly $\mu$-semistable rank-2 sheaves with Chern classes $c_1=0, \ c_2=m, \ c_3=2n$ (for more details, see \cite[Thm. 4.1, Prop. 4.2]{SRS}).

\begin{Lemma}\label{automorphisms}
For any sheaf $F$ from $\mathcal{V}(0,m,2n)$ we have that \emph{$h^{0}(F)=1$, $\text{dim End}(F) = 2$, $\text{Aut}(F) \simeq k^{*} \times k$}.
\end{Lemma}

\noindent\textit{Proof:} The extension (\ref{serre for non-stable}) immediately implies that $h^{0}(F)=1$. Also applying the functor $\text{Hom}(F, -)$ to (\ref{serre for non-stable}) we get
$$
0 \longrightarrow \text{Hom}(F, \mathcal{O}_{\mathbb{P}^{3}}) \longrightarrow \text{Hom}(F, F) \longrightarrow \text{Hom}(F, I_{Y}).
$$
So $\text{dim Hom}(F, F) \leq \text{dim Hom}(F, \mathcal{O}_{\mathbb{P}^{3}}) + \text{dim Hom}(F, I_{Y}) \leq 2 \ \text{dim Hom}(F, \mathcal{O}_{\mathbb{P}^{3}}) = 2 \ h^0(F) = 2$; second inequality holds because $I_{Y} \hookrightarrow \mathcal{O}_{\mathbb{P}^{3}}$ implies $\text{Hom}(F, I_{Y}) \hookrightarrow \text{Hom}(F, \mathcal{O}_{\mathbb{P}^{3}})$. On the other hand, the short exact sequence (\ref{serre for non-stable}) gives the following endomorphism 
\begin{equation}\label{endomorphism sigma}
\sigma: \ \ F \twoheadrightarrow I_{Y} \hookrightarrow \mathcal{O}_{\mathbb{P}^{3}} \hookrightarrow F,
\end{equation}
which is not a scalar multiplication. So, $\text{dim End}(F) \geq 2$, hence is 2. Therefore, the endomorphism algebra $\text{End}(F)$ has the following form
$$
\text{End}(F) \simeq \{ \lambda \text{Id} + \mu \sigma \ | \ \lambda, \mu \in k \}. 
$$
Finally, since $\sigma^{2} = 0$ for the corresponding automorphism group we have the isomorphism of groups
\begin{equation}\label{aut group}
\text{Aut}(F) = \{ \lambda \text{Id} + \mu \sigma \ | \ \lambda \in k^{*}, \ \mu \in k \} \simeq k^{*} \times k,
\end{equation} 
$$
\lambda \text{Id} + \mu \sigma \mapsto \Big( \frac{\mu}{\lambda}, \ \lambda \Big) \in k^{*} \times k. 
$$
\hfill $\Box$

\vspace{3mm}

\begin{Lemma}\label{ext_nonstable}
For $\mu$-semistable reflexive sheaves $F$ from $\mathcal{V}(0,m,2n)$ we have the following equalities
\begin{equation}\label{formula_for_ext1}
\emph{dim Ext}^{1}(F,F )= \emph{dim}~\mathcal{V}(0,m,2n),
\end{equation}
\begin{equation}\label{defect1}
\emph{dim Ext}^{2}(F,F) = h^{1}(N_{Y/\mathbb{P}^{3}}) + g.
\end{equation}
\end{Lemma}

\noindent\textit{Proof:} In order to show this apply the functor $\text{Hom}(-,F)$ to the triple (\ref{serre for non-stable}), then we obtain the exact sequence
\begin{equation}\label{seq F right}
0 \longrightarrow \text{Hom}(I_{Y},F) \longrightarrow \text{End}(F) \longrightarrow \text{H}^{0}(F) \longrightarrow
\end{equation}
$$
\longrightarrow \text{Ext}^{1}(I_{Y},F) \longrightarrow \text{Ext}^{1}(F,F) \longrightarrow \text{H}^{1}(F).
$$
From the fact that the triple (\ref{serre for non-stable}) is not splitting we can deduce that the canonical map $\text{Hom}(I_{Y}, \mathcal{O}_{\mathbb{P}^{3}}) \longrightarrow \text{Hom}(I_{Y}, F)$ is isomorphism, so we have
$$
\text{Hom}(I_{Y},F) \simeq \text{Hom}(I_{Y},\mathcal{O}_{\mathbb{P}^{3}}) \simeq k.
$$
Next, it is easy to see that $\text{H}^{0}(F) \simeq k$. It immediately implies that the morphism $\text{End}(F) \longrightarrow \text{H}^{0}(F)$ from the exact sequence (\ref{seq F right}) must be surjective. Since the curve $Y$ is smooth and irreducible we also have $\text{H}^{1}(F) \simeq \text{H}^{1}(I_{Y}) = 0$. Therefore, the exact sequence (\ref{seq F right}) implies the isomorphism
\begin{equation}\label{isom_for_ext}
\text{Ext}^{1}(F,F) \simeq \text{Ext}^{1}(I_{Y},F).
\end{equation}
Now after applying the functor $\text{Hom}(I_{Y},-)$ to the triple (\ref{serre for non-stable}) we have the exact sequence
\begin{equation}\label{seq_for_ext}
0 \longrightarrow \text{Hom}(I_{Y}, \mathcal{O}_{\mathbb{P}^{3}}) \overset{\simeq}{\longrightarrow} \text{Hom}(I_{Y}, F) \overset{0}{\longrightarrow} \text{Hom}(I_{Y}, I_{Y}) \longrightarrow
\end{equation}
$$
\longrightarrow \text{Ext}^{1}(I_{Y},\mathcal{O}_{\mathbb{P}^{3}}) \longrightarrow \text{Ext}^{1}(I_{Y},F) \longrightarrow \text{Ext}^{1}(I_{Y},I_{Y}) \longrightarrow \text{Ext}^{2}(I_{Y},\mathcal{O}_{\mathbb{P}^{3}}).
$$
Taking into account that $Y$ is smooth we conclude that $\text{Ext}^{2}(I_{Y},\mathcal{O}_{\mathbb{P}^{3}}) \simeq \text{Ext}^{3}(\mathcal{O}_{Y},\mathcal{O}_{\mathbb{P}^{3}}) \simeq \text{H}^{0}(\mathcal{O}_{Y}(-4)) = 0$ by Serre duality. Next, it is easy to check that $\text{Hom}(I_{Y},I_{Y}) \simeq k$ and
\begin{equation}\label{abc}
\text{Ext}^{1}(I_{Y}, \mathcal{O}_{\mathbb{P}^{3}}) \simeq \text{H}^{0}(\mathcal{E}xt^{1}(I_{Y}, \mathcal{O}_{\mathbb{P}^{3}})) \simeq
\end{equation}
$$
\simeq \text{H}^{0}(\mathcal{E}xt^{2}(\mathcal{O}_{Y},\mathcal{O}_{\mathbb{P}^{3}})) \simeq \text{H}^{0}(\omega_{Y}(4)).
$$
Moreover, we have $\text{Ext}^{1}(I_{Y},I_{Y}) \simeq \text{H}^{0}(\mathcal{E}xt^{1}(I_{Y},I_{Y})) \simeq \text{H}^{0}(N_{Y/\mathbb{P}^{3}})$. Substituting (\ref{isom_for_ext}), (\ref{abc}) into (\ref{seq_for_ext}) we obtain the exact sequence
\begin{equation}\label{seq_for_computation}
0 \longrightarrow k \longrightarrow \text{H}^{0}(\omega_{Y}(4)) \longrightarrow \text{Ext}^{1}(F,F) \longrightarrow \text{H}^{0}(N_{Y/\mathbb{P}^{3}}) \longrightarrow 0,
\end{equation}
which, together with (\ref{dim_strictmu}), immediately implies the equality (\ref{formula_for_ext1}). 

Since the curve $Y$ is smooth we have the formulas
$$
\chi(N_{Y/\mathbb{P}^{3}})=4 m, \ \ \ \chi(\omega_{Y}(4))=4 m + g - 1,
$$
$$
\text{deg}(\omega_{Y}(4))=4 m + 2g  - 2, \ \ \ h^{1}(\omega_{Y}(4))=0,
$$
so the equality (\ref{formula_for_ext1}) can be written in the following form
\begin{equation}\label{formula_for_ext1_new}
\text{dim Ext}^{1}(F,F)= h^{1}(N_{Y/\mathbb{P}^{3}}) + 8m + g - 2.
\end{equation}
Next, note that according to \cite[Prop. 3.4]{SRS} we have the following formula
\begin{equation}\label{Riemann-Roch}
\sum_{i=0}^{3} \text{dim Ext}^{i}(F,F)=-8 m + 4.
\end{equation}
Considering the remaining part of the exact sequence (\ref{seq F right}), namely,
$$
\text{Ext}^{3}(I_{Y},F) \longrightarrow \text{Ext}^{3}(F,F) \longrightarrow \text{H}^{3}(F)
$$
we can deduce that $\text{Ext}^{3}(F,F) = 0$ because of $\text{H}^{3}(F) \simeq \text{H}^{3}(I_{Y}) = 0$ and $\text{Ext}^{3}(I_{Y},F) \simeq \text{Ext}^{3}(I_{Y},I_{Y}) = 0$. On the other hand, from the Lemma \ref{automorphisms} we know that $\text{Hom}(F,F) \simeq k^{2}$. Taking into account (\ref{formula_for_ext1_new}), we obtain from (\ref{Riemann-Roch}) the formula (\ref{defect1}).
\hfill$\Box$

\vspace{3mm}

\noindent It is important that the proof of the irreducibility of the new components of $\mathcal{M}(k)$ which will be constructed in the next section using elementary transformations of reflexive sheaves $F$ is presented only for the case $\text{dim Ext}^{2}(F,F)=0$. By this reason and due to the formula (\ref{defect1}), we will consider only irreducible subschemes $\mathcal{V}_{m} \subset \mathcal{V}(0,m,4m-2)$ whose general points are sheaves obtained by Serre construction (\ref{serre for non-stable}) with smooth rational curves $Y$ of degree $m$ (it is easy to see that for such sheaves $F \in \mathcal{V}_{m}$ we have $c_3(F)=4m-2$). 

\begin{Remark}
\emph{For a sheaf $F \in \mathcal{V}_{m}$ the corresponding curve from the construction (\ref{serre for non-stable}) we will denote by $Y_{F}$. Also note that we have the inclusion $\text{Sing}(F) \subset Y_{F}$.}
\end{Remark}

\begin{Lemma}\label{triviality}
Let $F$ be an rank-$2$ $\mu$-semistable sheaf with $c_1(F) = 0$. Then for any $d \geq 1$, the restriction of $F$ to a general rational curve of degree $d$ in $\mathbb{P}^{3}$ is trivial.
\end{Lemma}

\noindent\textit{Proof:}
For $d = 1$, the assertion follows from the Grauert–M\"{u}lich Theorem \cite[Theorem 3.1.2]{HL}. For $d > 1$, we start by restriction to a general chain of m lines and then smooth out the chain of lines to a nonsingular rational curve of degree $d$. 

By a chain of lines we mean a curve $C_0 = l_1 \cup ... \cup l_d$ in $\mathbb{P}^{3}$ such that $l_1, ..., l_d$ are distinct lines and $l_i \cap l_j = \emptyset$ if and only if $|i-j| \leq 1$. It is well known (see e. g. \cite[Cor. 1.2]{HH}) that a chain of lines $C_0 = l_1 \cup ... \cup l_d$ in $\mathbb{P}^{3}$ considered as a reducible curve of degree $d$ can be deformed in a flat family with a smooth one-dimensional base $(\Delta, 0)$ to a nonsingular rational curve $C$. Making an \'{e}tale base change, we can obtain such a smoothing with a section.

By the case $d = 1$, the restriction of $F$ to a general line is trivial. By induction on $d$, we easily deduce that for a general chain of lines $C_0$, the restriction of $F$ to $C_0$ is also trivial: $F|_{C_{0}} \simeq \mathcal{O}^{\oplus  2}_{C_{0}}$, which is equivalent to saying that $F|_{l_i} \simeq \mathcal{O}^{\oplus 2}_{l_i}$ for all $i = 1, ..., d$. Choosing a smoothing $\{ C_{t} \}_{t \in \Delta}$ of $C_{0}$ with a section $t \mapsto x_{t} \in C_{t}$ as above, we remark that $F|_{C_{t}} \simeq \mathcal{O}_{C_{t}}(k_{t} ~ pt) \oplus \mathcal{O}_{C_{t}}( - k_{t} ~ pt)$ for some integer $k_{t}$ which may depend on $t$. The triviality of $F|_{C_{t}}$ is thus equivalent to the vanishing of $h^{0}(F|_{C_{t}}(-pt))$. Using the semi-continuity of $h^{0}(F|_{C_{t}}( - x_{t}))$, we see that $F|_{C_{t}}$ is trivial for general $t \in \Delta$.
\hfill$\Box$ 

\section{Construction of components}

Fix an arbitrary scheme $\mathcal{R}$ belonging to one of the families $\mathcal{S}_{a,b,c}$ or $\mathcal{V}_{m}$ described above. For simplicity of notation we will denote the Chern classes of a sheaf from $\mathcal{R}$ by $c_{i}(\mathcal{R}),~ i=1,2,3$. Similarly, fix some scheme $\mathcal{H}_{1}$ from the collection of the Hilbert schemes $\text{Hilb}_{d}, \ \text{Hilb}_{(d_1, d_2)}$, where the Hilbert scheme $\text{Hilb}_{d}$ parameterizes smooth irreducible rational curves of degree $d$ in $\mathbb{P}^{3}$ and $\text{Hilb}_{(d_1,d_2)}$ parametrizes smooth irreducible complete intersection curves of the form $S_{d_1} \cap S_{d_2}$, where $S_{d_1}, \ S_{d_2} \subset \mathbb{P}^{3}$ are surfaces of degree $d_1, \ d_2$, respectively. For the Hilbert schemes $\text{Hilb}_{(d_1,d_2)}$ we will assume that $d_{1} \leq d_{2}$ and $(d_{1}, d_{2}) \neq (1, 1), ~ (1, 2)$. The genus of curves from $\mathcal{H}_{1}$ we will denote by $g$ which is equal to zero for rational curves and $1+\frac{1}{2}d_1 d_2 (d_1 + d_2 - 4)$ for complete intersection curves. Next, denote by $\mathcal{H}_{0} = \text{Sym}^{s}_{*}(\mathbb{P}^{3})$ the open smooth subset of the Hilbert scheme parameterizing unions $W=\{x_1, ..., x_s \ | \ x_{i} \neq x_j \}$ of $s$ distinct points in $\mathbb{P}^{3}$. Also we impose the following restrictions on the choice of schemes $\mathcal{R}, \ \mathcal{H}_{1}$ and $\mathcal{H}_{0}$
\begin{equation}\label{cond on integers}
\left\{
\begin{array}{cl}
s < \frac{1}{2} c_3(\mathcal{R}) \ \ \ \text{if} \ \mathcal{H}_{1} = \text{Hilb}_{d}, \\
s \leq \frac{1}{2} c_3(\mathcal{R}) \ \ \ \text{if} \ \mathcal{H}_{1} = \text{Hilb}_{(d_1,d_2)}, \\
\mathcal{R}=\mathcal{V}_{m} \Rightarrow m<d.
\end{array}\right.
\end{equation}
The universal curves of the Hilbert schemes $\mathcal{H}_{0}$ and $\mathcal{H}_{1}$ we will denote by $\mathcal{Z}_{0} \subset \mathcal{H}_{0} \times \mathbb{P}^{3}$ and $\mathcal{Z}_{1} \subset \mathcal{H}_{1} \times \mathbb{P}^{3}$, respectively.

Since for a smooth projective curve invertible sheaves and rank-1 stable sheaves are the same objects, the relative Picard functor $\textbf{Pic}: (\textit{Sch}/\mathcal{H}_{1}) \longrightarrow (\textit{Sets})$ defined as
$$
\textbf{Pic}(T)=\{\text{$T$-flat invertible sheaves~} F \text{~on~} \mathcal{Z}_{1} \times_{\mathcal{H}_{1}} T \} / \text{Pic}(T)
$$
is equal to the relative Maruyama moduli functor for classifying stable sheaves which is corepresented by some $\mathcal{H}_{1}$-scheme (see \cite[Thm. 5.6]{Mar} or \cite[Thm. 4.3.7]{HL}). So the Picard functor $\textbf{Pic}$ is also corepresented by this $\mathcal{H}_{1}$-scheme which we denote by $\text{Pic}_{\mathcal{Z}_{1}/\mathcal{H}_1}$. Further we will only consider the component of the scheme $\text{Pic}_{\mathcal{Z}_{1}/\mathcal{H}_1}$ corresponding to the following Hilbert polynomial
$$P(k)=g-1+2d+n-s+dk.$$
We will denote this component just by $\mathcal{P}$. From the set-theoretical point of view the scheme $\mathcal{P}$ has the following form
$$
\mathcal{P}=\{ (C, L) \ | \ C \in \mathcal{H}_{1}, \ L \in \text{Pic}^{g-1+2 \text{deg}(C) +n-s}(C) \}.
$$
For the case of smooth rational curves we have the isomorphism $\mathcal{P} \simeq \mathcal{H}_{1}$ because $\text{Pic}^{g-1+2 \text{deg}(C) +n-s}(C)$ is trivial for any smooth rational curve $C$. Also it is obvious that the dimension of the scheme $\mathcal{P}$ can be computed by the formula
\begin{equation}\label{picard_dim}
\text{dim~}\mathcal{P} = \text{dim~} \mathcal{H}_{1} + \text{dim~ Jac}(C).
\end{equation}

\vspace{3mm}

\begin{Remark}
\emph{For simplicity of notation the sheaf $\mathcal{O}_{W} \oplus L$ for fixed elements $W \in \mathcal{H}_{0}, \ L \in \mathcal{P}$ we will denote just by $Q$ throughout the text. Also for arbitrary two sheaves $F$ and $Q$ we denote by $\text{Hom}_{e}(F,Q) \subset \text{Hom}(F,Q)$ the subset of surjective morphisms $F \twoheadrightarrow Q$ of $\text{Hom}(F,Q)$.} 
\end{Remark}

\vspace{3mm}

\begin{Lemma}
	The closed points of $\mathcal{R} \times \mathcal{P} \times \mathcal{H}_{0}$ satisfying the following conditions
	\begin{equation}\label{disjointness}
		C \cap W = \emptyset,
	\end{equation}
	\begin{equation}\label{supports}
	\mathcal{R} = \mathcal{S}_{a,b,c} \ \ \Rightarrow \ \ \emph{Sing}(F) \cap (C \sqcup W) = \emptyset,
	\end{equation}
	\begin{equation}\label{zeroes of section}
	\mathcal{R}=\mathcal{V}_{m} \ \ \Rightarrow \ \ Y_{F} \cap (C \sqcup W) = \emptyset,
	\end{equation}
	\begin{equation}\label{cond for h1}
		h^1(\mathcal{H}om(F,L))=0,
	\end{equation}
	\begin{equation}\label{cond epimorphism}
		\emph{Hom}_{e}(F, Q) \neq 0,
	\end{equation}
	\begin{equation}\label{cond defect}	
		h^{0}(\omega_{C}(4) \otimes L^{-2})=0
	\end{equation}
	form an open dense subset $\mathcal{B}$ of $\mathcal{R} \times \mathcal{P} \times \mathcal{H}_{0}$.
\end{Lemma}

\noindent\textit{Proof:} First of all, note that all these conditions are open, so we only need to prove that each of them is non-empty because of the irreducibility of the scheme $\mathcal{R} \times \mathcal{P} \times \mathcal{H}_{0}$. 

It is obvious that the conditions (\ref{disjointness})-(\ref{zeroes of section}) are non-empty because for a given sheaf $F \in \mathcal{R}$ we always can take the disjoint union $C \sqcup W$ away from $\text{Sing}(F)$ for the case $\mathcal{R}=\mathcal{S}_{a,b,c}$, or $Y_{F}$ for the case $\mathcal{R}=\mathcal{V}_{m}$. From Lemma \ref{triviality} it follows that the restriction of any sheaf from $\mathcal{R} = \mathcal{S}_{a,b,c}, \ \mathcal{V}_{m}$ on a general rational curve is trivial. Moreover, the first inequality of (\ref{cond on integers}) implies that $\text{deg}(L) > g - 1 + 2 ~ \text{deg}(C)>0$. Using these facts we can immediately conclude that the conditions (\ref{cond for h1})-(\ref{cond defect}) are non-empty for the case $\mathcal{H}_{1}=\text{Hilb}_{d}$.

Now let us prove that the conditions (\ref{cond for h1})-(\ref{cond defect}) are non-empty for the case $\mathcal{H}_{1}=\text{Hilb}_{(d_1,d_2)}(\mathbb{P}^{3})$ as well. In order to do this we consider a flat family $\textbf{C} :=\{ C_{t} \subset \mathbb{P}^{3}, \ t \in Y \} \subset Y \times \mathbb{P}^{3}$ of smooth curves parameterized by a smooth irreducible curve $Y$ with marked point $0 \in Y$ such that $C_{t} \in \mathcal{H}_{1}$ for $t \neq 0$, but $C_{0}=\bigcup\limits_{i=1}^{d_1 d_2} l_{i}$ is a union of $d_1 d_2$ projective lines $l_{1},...,l_{d_1 d_2}$ with at most double points. According to \cite[Lemma 20]{JMT2}, such family can be chosen with the property that there exists a sheaf $\widetilde{\textbf{L}}$ over $\textbf{C}$ satisfying the following
\begin{itemize}
	\item for $t \neq 0$: $\widetilde{L}_{t} := \widetilde{\textbf{L}}|_{t \times C_{t}} \in \text{Pic}^{g-1}(C_{t})$ and
    $$
	h^{0}(\widetilde{L}|_{C_{t}})=h^{1}(\widetilde{L}|_{C_{t}})=0, \ \ \ {\widetilde{L}|_{C_{t}}}^{\otimes 2} \neq \omega_{C_{t}};
	$$
	\item $\widetilde{L}_{0}=\bigoplus\limits_{i=1}^{d_1 d_2} \mathcal{O}_{l_{i}}(-1)$ is $\mu$-semistable with respect to the ample line bundle $\mathcal{O}_{\mathbb{P}^{3}}(1)|_{C_{0}}$.
\end{itemize}
Now fix a plane $H \subset \mathbb{P}^{3}$ which intersects $C_{0}$ at $d_1 d_2$ points, then $H$ transversally intersects the curve $C_{t}$ for any $t$ from some open subset $U \subset Y$ containing $0 \in Y$. Making an etale base change, we can assume that there is a section $x : U \longrightarrow \textbf{C}$ defined by $t \mapsto x_{t} \in C_{t} \cap H \subset C_{t}, \ t \in U$ such that $x_{0} \in l_{1}$. The section $x$ can be considered as the divisor $\{ x_{t} \}_{t \in U}$ on $\textbf{C}$ as well as $\{ H \cap C_{t} \}_{t \in U} \subset \textbf{C}$, so we can define the divisor $\textbf{D} := \{ H \cap C_{t} \}_{t \in U} + (d_1 d_2+n-s) \{ x_{t} \}_{t \in U}$ on $\textbf{C}$. Therefore, the sheaf $\textbf{L} := \widetilde{\textbf{L}}( \textbf{D})$ satisfies the following properties 
\begin{itemize}
	\item for $t \neq 0$: $L_{t} := \textbf{L}|_{t \times C_{t}} \in \text{Pic}^{g-1+2 d_1 d_2 +n-s}(C_{t})$
	\item $L_{0} = \mathcal{O}_{l_1}(d_1 d_2+n-s) \oplus \bigoplus\limits_{i=2}^{d_1 d_2} \mathcal{O}_{l_{i}}$.
\end{itemize}

Note that the restriction of a given sheaf $F$ from $\mathcal{R}$ on a general projective line is trivial due to Lemma \ref{triviality}. Conversely, using the induced action of the projective transformation group $\text{PGL}(4, k)$ on $\mathcal{R}$ we can state that restriction of a general sheaf from $\mathcal{R}$ on a given projective line is trivial. So there is a sheaf $[F] \in \mathcal{R}$ which is trivial on every line $l_{i}$ of the configuration $C_{0}$, i. e.
\begin{equation}\label{triv}
F|_{l_{i}} \simeq 2 \mathcal{O}_{l_i}, \ \ i=1,...,d_1 d_2.
\end{equation}
From this it is easy to see that
$$
h^{1}(\mathcal{H}om(F,{L}_{0}))=h^{1}(\mathcal{H}om(F,\mathcal{O}_{l_{1}}(d_1 d_2 + n-s)) \oplus \bigoplus\limits_{i=2}^{d_1 d_2} h^{1}(\mathcal{H}om(F,\mathcal{O}_{l_{i}}))=0.
$$
Taking into account the upper-semicontinuity we can conclude that this equality holds for $L_{t}$, where $t$ belongs to some open subset of $U$. Since $([F], L_{t}) \in \mathcal{R} \times \mathcal{P}$ for $t \neq 0$ this proves that the property (\ref{cond for h1}) is non-empty.

Next, since the sheaf $F$ is locally-free free along the support of the sheaf $L_{t}$ and as it was proved the equality $h^{i \geq 1}(\mathcal{H}om(F,L_{0})) = 0$ holds, any epimorphism $F \twoheadrightarrow L_{0}$ can be extended to epimorphism $F \twoheadrightarrow L_{t}$ (see \cite[Lemma 7.1]{JMT1}). So we have $\text{Hom}_{e}(F, L_{t}) \neq 0$, and, obviously, $\text{Hom}_{e}(F,\mathcal{O}_{W} \oplus L_{t}) \neq 0$ for any $W \in \mathcal{H}_{0}$ not intersecting $C_{t} \sqcup \text{Sing}(F)$.

Finally, let us prove that the condition (\ref{cond defect}) is non-empty. Note that for any pair $(C,L) \in \mathcal{P}$ we have the following equality
$$
\text{deg}(\omega_{C}(4) \otimes L^{-2})=2g-2+4 d_1 d_2-2(g-1+2 d_1 d_2 +n-s)=2(s-n).
$$
The second inequality of (\ref{cond on integers}) means that $s \leq n$. So if $s<n$ then the condition (\ref{cond defect}) is obviously satisfied. On the other hand, if $s=n$ then the line bundle $L$ can be chosen in such way that $L \simeq \widetilde{L}(2)$, where $\widetilde{L}$ is not a theta-characteristic, i. e. $\widetilde{L}^{\otimes 2} \neq \omega_{C}$. Hence $\omega_{C}(4) \otimes L^{-2}$ is non-trivial line bundle of degree 0 so it has no nonzero global sections.
\hfill$\Box$

\vspace{3mm}

\begin{Lemma}\label{equivalence}
Two different triples $(F, W \sqcup C, L), \ (\widetilde{F}, \widetilde{W} \sqcup \widetilde{C}, \widetilde{L}) \in \mathcal{B}$ give the same isomorphism class $[E \overset{\phi}{\simeq} \widetilde{E}] \in \mathcal{M}(m+d)$ if and only if there exist isomorphisms $\psi \in \emph{Hom}(F,\widetilde{F}), \ \zeta \in \emph{Hom}(Q,\widetilde{Q})$ which complete the commutative diagram
\begin{equation}\label{completion of diagram}
\begin{tikzcd}
0 \rar & E \rar{\xi} \dar{\phi} & F \rar \dar{\psi} & Q \rar \dar{\zeta} & 0 & \\
0 \rar & \widetilde{E} \rar{\widetilde{\xi}} & \widetilde{F} \rar & \widetilde{Q} \rar & 0 &
\end{tikzcd}
\end{equation}
\end{Lemma}
\noindent\textit{Proof}: See \cite[Cor. 1.5]{SRS}.
\hfill$\Box$

\vspace{3mm}

\begin{Lemma}\label{stability}
For any triple $(F,W \sqcup C, L) \in \mathcal{B}$ and surjective morphism $\phi \in \emph{Hom}_{e}(F, Q)$ the kernel sheaf $E:=\emph{ker}~\phi$ is stable.
\end{Lemma}
\noindent\textit{Proof}: It is obvious that the sheaf $E$ is $\mu$-semistable. Moreover, it has no torsion and $c_{1}(E)=0$, so in order to prove its stability we can consider only subsheaves $G \subset E$ which are sheaves of ideals $I_{\Delta}$ of some subschemes $\Delta \subset \mathbb{P}^{3}, \ \text{dim}~\Delta \leq 1$. Since taking double dual sheaf is functorial we have the following commutative diagram
\begin{equation}
\begin{tikzcd}[column sep=small]
& 0 \dar & 0 \dar \\
0 \rar & I_{\Delta} \rar \dar & \mathcal{O}_{\mathbb{P}^{3}} \dar \\
0 \rar & E \rar & E^{\vee\vee} \\
\end{tikzcd}
\end{equation}
which implies that $h^{0}(E^{\vee\vee}) > 0$. On the other hand, as it was shown in the proof of the previous lemma there is the isomorphism $E^{\vee \vee} \simeq F$. However, for the case $\mathcal{R}=\mathcal{S}_{a,b,c}$ the corresponding sheaf $[F] \in \mathcal{R}$ has no nonzero global sections because it would contradict its stability. Therefore, for this case the sheaf $E$ is stable.

Next, consider the case $\mathcal{R}=\mathcal{V}_{m}$. Note that from the construction (\ref{serre for non-stable}) it follows that $h^0(F)=1$. So the diagram above can be written in the following form
\begin{equation}
\begin{tikzcd}[column sep=small]
& 0 \dar & 0 \dar & 0 \dar \\
0 \rar & I_{\Delta} \rar \dar & \mathcal{O}_{\mathbb{P}^{3}} \rar \dar & \mathcal{O}_{\Delta} \rar \dar & 0 \\
0 \rar & E \rar \dar & F \rar \dar & \mathcal{O}_{W} \oplus L \rar & 0 \\
0 \rar & T \rar \dar & I_{Y} \dar & & \\
& 0 & 0 &
\end{tikzcd}
\end{equation}
From this we immediately conclude that $\Delta \subset W \sqcup C$. Note that $C$ is an irreducible curve and $L$ is a locally-free $\mathcal{O}_{C}$-sheaf, so it cannot have $0$-dimensional subsheaf. Therefore, in the case $\text{dim} \ \Delta = 0$, the composition $\mathcal{O}_{\Delta} \hookrightarrow \mathcal{O}_{W} \oplus L \overset{pr_{2}}{\twoheadrightarrow} L$ is zero and $\Delta \cap C \neq \emptyset$ is impossible. Hence, only the following cases are possible: $\Delta \cap C=\emptyset$ or $C$. The first case leads to contradiction because $I_{Y}|_{C} \simeq \mathcal{O}_{C}$ and $\text{deg}~L>0$, so there is no surjective morphism $I_{Y} \twoheadrightarrow L$. The second case is not destabilizing due to the third inequality of (\ref{cond on integers}) and the following formula
$$
\frac{1}{2}P(E)-P(I_{C \sqcup W'}) =  \frac{\text{deg}(C)-m}{2}k + \text{const}.
$$
Therefore, for the case $\mathcal{R}=\mathcal{V}_{m}$ the sheaf $E$ is also stable.
\hfill$\Box$

\vspace{3mm}

From Lemma \ref{equivalence} it follows that an epimorphism $\phi \in \text{Hom}_{e}\big( F, Q \big)$ defines the isomorphism class $[\text{ker}~\phi]$ up to natural action of $\text{Aut}(F) \times \text{Aut}(Q)$ on $\text{Hom}_{e}\big( F, Q \big)$, i. e. $(\psi, \zeta) \phi = \zeta \circ \phi \circ \psi^{-1}$ for $(\psi, \zeta) \in \text{Aut}(F) \times \text{Aut}(Q)$. In other words, the element $[\phi]$ of the orbit space $\text{Hom}_{e}\big( F, Q \big) / \Big( \text{Aut}(F) \times \text{Aut}(Q) \Big)$, which we will consider as a set, uniquely defines the isomorphism class $[\text{ker}~\phi]$. This fact, together with Lemma \ref{stability}, implies that the elements of the following set of data of elementary transformations
\begin{equation}
\mathcal{Q}:=\Big\{([F], C \sqcup W, L, [\phi]) \ | \ ([F], C \sqcup W, L) \in \mathcal{B},
\end{equation}
$$
\ [\phi] \in \text{Hom}_{e}\big( F, Q \big) / \Big( \text{Aut}(F) \times \text{Aut}(Q) \Big) \Big\}
$$
\noindent are in one-to-one correspondence with some subset of closed points of the moduli scheme $\mathcal{M}(m+d)$. Note that the vector space $\text{Hom}\big( F, Q \big)$ has the following direct decomposition
\begin{equation}\label{decomposition}
\text{Hom}\big( F, Q \big)=\text{Hom}\big( F, L \big) \oplus \text{Hom}\big( F, \mathcal{O}_{x_1} \big) \oplus ... \oplus \text{Hom}\big( F, \mathcal{O}_{x_s} \big).
\end{equation}
Moreover, all non-trivial morphisms $F \longrightarrow \mathcal{O}_{x_i}$ are surjective, so we have 
$$
\text{Hom}_{e}\big( F, \mathcal{O}_{x_i} \big) = \text{Hom}\big( F, \mathcal{O}_{x_i} \big) \setminus 0 = k^{2} \setminus 0, \ \ \ i=1,...,s.
$$
Next, since the sheaves $L, \mathcal{O}_{x_1}, ..., \mathcal{O}_{x_s}$ are simple and their supports are disjoint we have the isomorphism 
$$
\text{Aut}(Q) \simeq  \text{Aut}(L) \times \text{Aut}(\mathcal{O}_{x_{1}}) \times ... \times \text{Aut}(\mathcal{O}_{x_{s}}) \simeq (k^{*})^{s+1} 
$$
which obviously respects the decomposition (\ref{decomposition}), so we have the following equality
\begin{equation}\label{proj decomposition}
\text{Hom}_{e}\big( F, Q \big) / \text{Aut}(Q) \simeq \mathbb{P}\text{Hom}_{e}\big( F, L \big) \times \prod_{i=1}^{s} \mathbb{P} \text{Hom}_{e}\big( F, \mathcal{O}_{x_i} \big) \simeq
\end{equation}
$$
\simeq \mathbb{P}\text{Hom}_{e}\big( F, L \big) \times (\mathbb{P}^{1})^{\times s} \overset{open}{\lhook\joinrel\longrightarrow} \mathbb{P}\text{Hom}\big( F, L \big) \times (\mathbb{P}^{1})^{\times s}. 
$$
From this it follows that
\begin{equation}
\text{Hom}_{e}\big( F, Q \big) / \Big( \text{Aut}(F) \times \text{Aut}(Q) \Big) = 
\end{equation}
$$
=\Big( \mathbb{P}\text{Hom}_{e}\big( F, L \big) \times \prod_{i=1}^{s} \mathbb{P} \text{Hom}_{e}\big( F, \mathcal{O}_{x_i} \big) \Big) / \text{PAut}(F),
$$
\noindent where $\text{PAut}(F) :=  \text{Aut}(F) / \{ \lambda \cdot \text{Id}, \ \lambda \in k^{*} \}$ is the quotient group by homotheties. 

For the case $\mathcal{R} = \mathcal{S}_{a, b, c}$ the automorphism group $\text{Aut}(F)$ is generated by the homotheties, so the group $\text{P}\text{Aut}(F)$ is trivial. On the other hand, for the case $\mathcal{R} = \mathcal{V}_{m}$ the automorphism group of the sheaf $F$ is of the form (\ref{aut group}), so the group $\text{PAut}(F) \simeq k$ is not trivial. Let us show how the automorphism group $\text{Aut}(F)$ acts on the vector space $\text{Hom}(F,Q)$. From the exact triple (\ref{serre for non-stable}), the conditions (\ref{zeroes of section}) and (\ref{cond for h1}) we have the following exact triple
$$
0 \longrightarrow \text{Hom}(I_{Y},Q) \longrightarrow \text{Hom}(F,Q) \longrightarrow  \text{Hom}(\mathcal{O}_{\mathbb{P}^{3}},Q) \longrightarrow 0.
$$
Note that $\text{Hom}(I_{Y}, Q) \simeq \text{Hom}(\mathcal{O}_{\mathbb{P}^{3}}, Q)=:V$, so we have the isomorphism $\text{Hom}(F,Q) \simeq V \oplus V$. It is easy to see that the endomorphism $\sigma \in \text{End}(F)$ induces the following action on $V \oplus V$ by sending $(x, y) \in V \oplus V$ to $(y, 0)$. Moreover, since there are no surjections $I_{Y} \twoheadrightarrow Q$ we have that $\text{Hom}_{e}(F,Q) \cap \text{ker}~\sigma = \emptyset$. In particular, it means that the induced action of $\text{PAut}(F)$ on $\text{Hom}_{e}\big( F, Q \big) / \text{Aut}(Q)$ is free.

\begin{Lemma}
There exists an irreducible closed subscheme $\overline{\mathcal{C}}$ of $\mathcal{M}(m+d)$ and a dense subset $\mathcal{C} \subset \overline{\mathcal{C}}$ whose closed points are in one-to-one correspondence with the elements of the set $\mathcal{Q}$. The dimension of $\overline{\mathcal{C}}$ can be computed by the following formula
\begin{equation}\label{dim}
\emph{dim}~\overline{\mathcal{C}} = \emph{dim}~\mathcal{R} + \emph{dim~} \mathcal{H}_{0} + \emph{dim~} \mathcal{P} +
\end{equation}
$$
+\emph{dim~}\emph{Hom}(F, Q) / \emph{Aut}(Q) - \emph{dim }\emph{PAut}(F).
$$
\end{Lemma}

\noindent\textit{Proof:} Since the Hilbert scheme $\mathcal{H}_{1}$ parameterizes only smooth connected curves, there exists a Poincar\'{e} sheaf $\textbf{L}$ on $\mathcal{P} \times_{\mathcal{H}_{1}} \mathcal{Z}_{1}$. Next, for the case $\mathcal{R} = \mathcal{S}_{a,b,c}$ there exists an etale surjective morphism $\xi : \widetilde{\mathcal{R}} \longrightarrow \mathcal{R}$ and a sheaf $\textbf{F}$ on $\widetilde{\mathcal{R}} \times \mathbb{P}^{3}$ such that $\textbf{F}|_{t \times \mathbb{P}^{3}} \simeq F_{\xi(t)}$, where $[F_{\xi(t)}]$ is the isomorphism class of the sheaf defined by the point $\xi(t) \in \mathcal{R}$. The etale morphism $\xi$ can be obtained in the following way. Recall the construction of the moduli space $\mathcal{R}$ (see \cite[Thm 4.3.7]{HL}). Namely, $\mathcal{R}$ is obtained as a GIT-quotient $p : \frak{Q} \longrightarrow \frak{Q} //  GL(N) = \mathcal{R}$ for an appropriately chosen open subset $\frak{Q}$ of the Quot-scheme $\text{Quot}_{\mathbb{P}^{3}}(\mathcal{O}_{\mathbb{P}^{3}}(-m)^{\oplus N},P)$, where $P$ is the corresponding Hilbert polynomial, $N=P(m)$ and $m$ large enough. Since the sheaves from $\mathcal{R}$ are stable, the quotient $p$ is a principal $GL(N)$-bundle (see \cite[Cor. 4.3.5]{HL}). By the definition this means that there exists an etale surjective morphism $\xi : \widetilde{\mathcal{R}} \twoheadrightarrow \mathcal{R}$ such that $\frak{Q} \times_{\mathcal{R}} \widetilde{\mathcal{R}}$ is isomorphic to the direct product $\widetilde{\mathcal{R}} \times GL(N)$. On the other hand, it is well-known that there exists the universal sheaf $\mathcal{F}$ over $\frak{Q} \times \mathbb{P}^{3}$. Denote by $\widetilde{\mathcal{F}}$ the corresponding pullback of $\mathcal{F}$ to $( \frak{Q} \times_{\mathcal{R}} \widetilde{\mathcal{R}} ) \times \mathbb{P}^{3}$. Next, let $\textbf{F}$ be the restriction of $\widetilde{\mathcal{F}}$ on $\widetilde{\mathcal{R}} \times \mathbb{P}^{3} \hookrightarrow (\widetilde{\mathcal{R}} \times GL(N) ) \times \mathbb{P}^{3} \simeq (\frak{Q} \times_{\mathcal{R}} \widetilde{\mathcal{R}}) \times \mathbb{P}^{3}$, where the first inclusion is given by fixing arbitrary point in $GL(N)$. It is obvious that all restrictions $\textbf{F}|_{t_{i} \times \mathbb{P}^{3}}, \ t_{i} \in \xi^{-1}(t_0)$ of the sheaf $\textbf{F}$ are isomorphic. Also note that since the scheme $\mathcal{R}$ is smooth, by taking the irreducible component of $\widetilde{\mathcal{R}}$, we can assume that $\widetilde{\mathcal{R}}$ is smooth and irreducible, and it covers some open dense subset of the scheme $\mathcal{R}$. For the case $\mathcal{R} = \mathcal{V}_{m}$ we can assume that $\widetilde{\mathcal{R}} = \mathcal{R}$ because there exists the universal sheaf $\mathbf{F}$ on $\mathcal{R} \times \mathbb{P}^{3}$ (see \cite{Lange}).  

From \cite[Lemma 4.5]{Str} once can deduce that the scheme $\textbf{P}(\textbf{F})$ is irreducible and reduced. The symmetric group $G = S_{s}$ acts on $\prod_{i=1}^{s} \textbf{P}(\textbf{F})$ by permutations of factors, and the $s$-fold fibered product $\textbf{P}(\textbf{F}) \times_{\widetilde{\mathcal{R}}} \cdot\cdot\cdot \times_{\widetilde{\mathcal{R}}} \textbf{P}(\textbf{F})$ naturally embeds in $\prod_{i=1}^{s} \textbf{P}(\textbf{F})$ as a $G$-invariant subscheme. Now consider the following integral scheme
$$
\text{Sym}^{s}_{\widetilde{\mathcal{R}}}(\textbf{P}(\textbf{F})) := \Big( \textbf{P}(\textbf{F}) \times_{\widetilde{\mathcal{R}}} \cdot\cdot\cdot \times_{\widetilde{\mathcal{R}}} \textbf{P}(\textbf{F}) \Big) / G.
$$
Since there is the projection $\textbf{P}(\textbf{F}) \longrightarrow \widetilde{\mathcal{R}} \times \mathbb{P}^{3}$, we also have the two natural projections 
$$
\text{Sym}^{s}_{\widetilde{\mathcal{R}}}(\textbf{P}(\textbf{F})) \longrightarrow \widetilde{\mathcal{R}}, \ \ \ \text{Sym}^{s}_{\widetilde{\mathcal{R}}}(\textbf{P}(\textbf{F})) \longrightarrow \text{Sym}^{s}(\mathbb{P}^{3}), 
$$ 
so we can define the following surjective morphism
$$
\text{Sym}^{s}_{\widetilde{\mathcal{R}}}(\textbf{P}(\textbf{F})) \longrightarrow \widetilde{\mathcal{R}} \times \text{Sym}^{s}(\mathbb{P}^{3}).
$$
Therefore, we can consider the fiber product of the following form
$$
Y:=\text{Sym}^{s}_{\widetilde{\mathcal{R}}}(\textbf{P}(\textbf{F})) \times_{\widetilde{\mathcal{R}} \times \text{Sym}^{s}(\mathbb{P}^{3})} (\mathcal{B} \times_{\mathcal{R}} \widetilde{\mathcal{R}}).
$$
Next, define the sheaf $\tau:=p_{*} \mathcal{H}om(q_{1}^{*}\textbf{F}, \ q_{2}^{*} \textbf{L})$ over $Y$, where $p, \ q_1, \ q_2$ are the natural projections included into the following diagram 
\begin{center}
$$\begin{CD}
Y   @<{p}<< Y \times_{\mathcal{H}_{1}} \mathcal{Z}_{1}  @>{q_1}>>   \widetilde{\mathcal{R}} \times \mathbb{P}^{3} \\
    @.      @V{q_2}VV \\
    @.      \mathcal{P} \times_{\mathcal{H}_{1}} \mathcal{Z}_{1}
\end{CD}$$
\end{center}
Assume that a point $y \in Y$ projects to the triple $(F, W \sqcup C, L) \in \mathcal{B}$, then the fiber $\tau_{y} \otimes k(y)$ of the sheaf $\tau$ over the point $y$ is isomorphic to $\text{Hom}(F,L)$. Due to the condition (\ref{cond for h1}) we have the equality $\chi(\mathcal{H}om(F,L)) = h^0(\mathcal{H}om(F,L))$. On the other hand, if $\mathcal{R}=\mathcal{S}_{a,b,c}$  then from the condition (\ref{supports}) and the exact triple (\ref{reflexive series}) it follows that $\chi(\mathcal{H}om(F,L))$ depends only on the Euler characteristics of the line bundle $L$. More precisely, we have that
$$
\chi(\mathcal{H}om(F,L)) = (a+b+c+2) \cdot \chi(L(k)) -
$$
$$
 - a \cdot \chi(L(k+3)) - b \cdot \chi(L(k+2)) - c \cdot \chi(L(k+1)).
$$
Since $\chi(L(k)), k \in \mathbb{Z}$ is constant for all $L \in \mathcal{P}$ we can conclude that all fibers $\tau_{y} \otimes k(y), \ y \in Y$ are of the same dimension. Similarly, for the case $\mathcal{R}=\mathcal{V}_{m}$ due to the condition (\ref{zeroes of section}) and the triple (\ref{serre for non-stable}) we can obtain the following formula
$$
\chi(\mathcal{H}om(F,L)) = \chi(\mathcal{H}om(\mathcal{O}_{\mathbb{P}^{3}}, L)) + \chi(\mathcal{H}om(I_{Y_{F}}, L)) = 2 \cdot \chi(L)
$$
which as previously implies that the fibers $\tau_{y} \otimes k(y), \ y \in Y$ are of the same dimension. From the construction it follows that the scheme $Y$ is reduced, so the sheaf $\tau$ is actually locally-free. Therefore, it can be viewed as a vector bundle.

Now consider the projective bundle $\textbf{P}(\tau^{\vee})$ associated to the vector bundle $\tau$. If the point $u \in \mathcal{B} \times_{\mathcal{R}} \widetilde{\mathcal{R}}$ projects to the point $(F,W \sqcup C, L) \in \mathcal{B}$, then the fiber of the projection $\textbf{P}(\tau^{\vee}) \longrightarrow \mathcal{B} \times_{\mathcal{R}} \widetilde{\mathcal{R}}$ over the point $u$ is the direct product of projective spaces $\mathbb{P}\text{Hom}(F,L) \times \prod_{i=1}^{s} \mathbb{P}\text{Hom}(F,\mathcal{O}_{x_i})$ which is isomorphic to $\text{Hom}(F,Q) / \text{Aut}(Q)$ according to (\ref{proj decomposition}). From the construction it follows that the dimension of $\textbf{P}(\tau^{\vee})$ can be computed by the following formula 
\begin{equation}\label{dim of proj}
\text{dim}~\textbf{P}(\tau^{\vee}) = \text{dim}~\mathcal{R} + \text{dim~} \mathcal{H}_{0} + \text{dim}~ \mathcal{P} + \text{dim Hom}(F, Q) / \text{Aut}(Q).
\end{equation}
Let $\frak{E} \subset \textbf{P}(\tau^{\vee})$ be the open dense subset of $\textbf{P}(\tau^{\vee})$ consisted from the classes of surjective morphisms $[F \twoheadrightarrow Q]$. Any point $q \in \frak{E}$ determines the isomorphism class of the sheaf $[E_{q}]:=[\text{ker} \ \psi_{q}]$, where $[\psi_{q}] \in \text{Hom}_{e}(F,Q) / \text{Aut}(Q)$. As in \cite[Prop. 6.4 ]{JMT1}, one can show that the family $\{ E_{q}, \ q \in \frak{E} \}$ globalizes in a standard way to the universal sheaf $\textbf{E}$ over $\frak{E} \times \mathbb{P}^{3}$. Next, by the construction and by the definition of moduli scheme, the sheaf $\textbf{E}$ defines the modular morphism $\Phi: \frak{E} \longrightarrow \mathcal{M}(m+d), \ q \mapsto [E_{q}=\text{ker} \ \psi_{q}]$. Now consider the image $\mathcal{C}:=\text{im}(\Phi)$ of the morphism $\Phi$ and its scheme-theoretic closure $\overline{\mathcal{C}} \subset \mathcal{M}(m+d)$. Note that the scheme $\frak{E}$ is irreducible, so the scheme $\overline{\mathcal{C}}$ is also irreducible. Moreover, the morphism $\Phi$ is flat over an dense open subset of $\mathcal{C}$. In particular, it means that for the general point $[E] = \Phi(x), \ x \in \frak{E}$ we have the following formula for the dimension
\begin{equation}\label{dim of image}
\text{dim}_{[E]}~ \overline{\mathcal{C}} = \text{dim}_{x}~\frak{E} - \text{dim}_{x}~\Phi^{-1}([E]).
\end{equation}
From the Lemma \ref{equivalence} it follows that $\Phi(x)=\Phi(y)$ if and only if the points $x, y \in \textbf{P}(\tau^{\vee})$ projects to the same tuple $(F,W \sqcup C, L) \in \mathcal{B}$ and the corresponding equivalence classes $[\phi_{x}], \ [\phi_{y}] \in \text{Hom}(F,Q) / \text{Aut}(Q)$ differs by the action of the group $\text{P}\text{Aut}(F)$ which is free. Therefore, the fiber of $\Phi^{-1}([E])$ is isomorphic to the disjoint union of the finite number of copies of the group $\text{P}\text{Aut}(F), \ F \simeq E^{\vee \vee}$. It implies that the set of closed points of $\mathcal{C}$ is isomorphic to $\mathcal{Q}$. Moreover, from the formulas (\ref{dim of proj}) and (\ref{dim of image}) it follows that the dimension of the scheme $\overline{\mathcal{C}} \subset \mathcal{M}(m+d)$ can be computed by the formula (\ref{dim}).
\hfill$\Box$

\vspace{3mm}

\section{Irreducibility of components}

\begin{Theorem}\label{Main result}
For general sheaf $[E]$ of the closed subscheme $\overline{\mathcal{C}} \subset \mathcal{M}(m+d)$ we have the equality
$$
\emph{dim}~T_{[E]}\mathcal{M}(m+d)=\emph{dim}~\overline{\mathcal{C}}.
$$
Therefore, the subscheme $\overline{\mathcal{C}}$ is an irreducible component of the moduli scheme $\mathcal{M}(m+d)$.
\end{Theorem}

\noindent\textit{Proof:} For the computation of the dimension of the tangent space of the moduli scheme $\mathcal{M}(m+d)$ at the point $[E]$ defined above, we use the standard fact of deformation theory, $T_{[E]}\mathcal{M}(m+d) \simeq \text{Ext}^{1}(E,E)$ for a stable sheaf $E$, and the local-to-global spectral sequence $\text{H}^p(\mathcal{E}xt^{q}(G,E)) \Rightarrow \text{Ext}^{p+q}(G,E)$ for any sheaf $G$, which yields the following exact sequence
\begin{equation}
0 \longrightarrow \text{H}^{1}(\mathcal{H}om(G,E)) \longrightarrow \text{Ext}^{1}(G,E) \longrightarrow \text{H}^{0}(\mathcal{E}xt^{1}(G,E)) \overset{\phi}{\longrightarrow}
\end{equation}
$$
\overset{\phi}{\longrightarrow} \text{H}^{2}(\mathcal{H}om(G,E)) \longrightarrow \text{Ext}^{2}(G,E).
$$

According to our construction, the general sheaf $[E] \in \overline{\mathcal{C}}$ fits into the exact triple of the following form
\begin{equation}\label{main}
0 \longrightarrow E \longrightarrow F \longrightarrow Q \longrightarrow 0.
\end{equation}
Note again that for general sheaves from $\mathcal{R} \in \{\mathcal{S}_{a,b,c}, \ \mathcal{V}_{m}\}$ we have the equality $\text{Ext}^{2}(F,F) = 0$ (see \cite[Lemma 5]{JMT2}). Moreover, from (\ref{cond for h1}) and (\ref{supports}) it follows that $\text{Ext}^{1}(F,Q) = 0$, so the following triple
\begin{equation}
\text{Ext}^{1}(F,Q) \longrightarrow \text{Ext}^{2}(F,E) \longrightarrow \text{Ext}^{2}(F,F)
\end{equation}
\noindent yields $\text{Ext}^{2}(F,E) = 0$. Taking into account that the sheaf $\mathcal{H}om(E,E)/\mathcal{H}om(F,E) \hookrightarrow \mathcal{E}xt^{1}(Q, E)$ has the dimension at most $1$, we have that the map $\text{H}^{2}(\mathcal{H}om(F, E)) \longrightarrow \text{H}^{2}(\mathcal{H}om(E, E))$ is surjective. Therefore, we obtain the commutative diagram
\begin{equation*}
\begin{tikzcd}[column sep=small]
\text{H}^{0}(\mathcal{E}xt^{1}(F, E)) \rar[two heads] \dar & \text{H}^{2}(\mathcal{H}om(F, E)) \rar \dar[two heads] & 0 \dar \\
\text{H}^{0}(\mathcal{E}xt^{1}(E, E)) \rar{\phi} & \text{H}^{2}(\mathcal{H}om(E, E)) \rar & \text{Ext}^{2}(E,E)
\end{tikzcd}
\end{equation*}
\noindent from which it follows that the morphism $\phi$ is surjective. Consequently, we have the following formula
\begin{equation}\label{equality for ext1}
\text{dim~}\text{Ext}^{1}(E,E)=h^{0}(\mathcal{E}xt^{1}(E,E))+h^{1}(\mathcal{H}om(E,E))-h^{2}(\mathcal{H}om(E,E)),
\end{equation}
\noindent and an analogous formula for the sheaf $F$
\begin{equation}\label{equality for ext1f}
\text{dim~}\text{Ext}^{1}(F,F)=h^{0}(\mathcal{E}xt^{1}(F,F))+h^{1}(\mathcal{H}om(F,F))-h^{2}(\mathcal{H}om(F,F)).
\end{equation}

Applying the functor $\mathcal{H}om(-,E)$ to the triple (\ref{main}) we obtain the following exact sequence
\begin{equation}
0 \longrightarrow \mathcal{H}om(Q,E) \longrightarrow \mathcal{H}om(F,E)  \longrightarrow \mathcal{H}om(E,E) \longrightarrow
\end{equation}
$$
\longrightarrow \mathcal{E}xt^{1}(Q,E) \overset{0}{\longrightarrow} \mathcal{E}xt^{1}(F,E) \longrightarrow \mathcal{E}xt^{1}(E,E) \longrightarrow
$$
$$
\longrightarrow \mathcal{E}xt^{2}(Q,E) \overset{0}{\longrightarrow} \mathcal{E}xt^2(F,E).
$$
Since $E$ is torsion-free sheaf we have $ \mathcal{H}om(Q,E)=0$. Since the sheaf $\mathcal{E}xt^{i \geq 1}(F,E)$ is supported on the subset $\text{Sing}(F)$, the condition (\ref{supports}) implies that the sheaves $\mathcal{E}xt^{1,2}(Q,E)$ and $\mathcal{E}xt^{1,2}(F,E)$ have disjoint supports, so the morphisms $\mathcal{E}xt^{1}(Q,E) \longrightarrow \mathcal{E}xt^{1}(F,E)$ and $\mathcal{E}xt^{2}(Q,E) \longrightarrow \mathcal{E}xt^2(F,E)$ are equal to zero. Therefore we obtain the following triples
\begin{equation}\label{triple2}
0 \longrightarrow \mathcal{H}om(F,E)  \longrightarrow \mathcal{H}om(E,E) \longrightarrow \mathcal{E}xt^{1}(Q,E) \longrightarrow 0,
\end{equation}
\begin{equation*}
0 \longrightarrow \mathcal{E}xt^{1}(F,E) \longrightarrow \mathcal{E}xt^{1}(E,E) \longrightarrow \mathcal{E}xt^{2}(Q,E) \longrightarrow 0.
\end{equation*}
For the same reason (\ref{supports}) implies that the last triple splits, so we have the isomorphism
\begin{equation}\label{isom1 for ext1}
\mathcal{E}xt^{1}(E,E) \simeq \mathcal{E}xt^{1}(F,E) \oplus \mathcal{E}xt^{2}(Q,E).
\end{equation}

Now apply the functor $\mathcal{H}om(F,-)$ to the triple (\ref{main})
\begin{equation}
0 \longrightarrow \mathcal{H}om(F,E) \longrightarrow \mathcal{H}om(F,F) \longrightarrow \mathcal{H}om(F,Q) \overset{0}{\longrightarrow}
\end{equation}
$$
\overset{0}{\longrightarrow} \mathcal{E}xt^{1}(F,E) \longrightarrow \mathcal{E}xt^{1}(F,F) \longrightarrow \mathcal{E}xt^{1}(F,Q).
$$
Again from (\ref{supports}) it follows that $\mathcal{E}xt^{1}(F,Q) = 0$ and the supports of the sheaves $\mathcal{H}om(F,Q), \ \mathcal{E}xt^{1}(F,E)$ are disjoint, so the morphism $\mathcal{H}om(F,Q) \longrightarrow \mathcal{E}xt^{1}(F,E)$ is equal to zero. Therefore, we obtain the following exact triple and isomorphism
\begin{equation}\label{triple1}
0 \longrightarrow \mathcal{H}om(F,E) \longrightarrow \mathcal{H}om(F,F) \longrightarrow \mathcal{H}om(F,Q) \longrightarrow 0,
\end{equation}
\begin{equation}\label{isom2 for ext1}
 \mathcal{E}xt^{1}(F,E) \simeq \mathcal{E}xt^{1}(F,F).
\end{equation}

Next, apply the functor $\mathcal{H}om(Q,-)$ to the triple (\ref{main})
\begin{equation}
0 \longrightarrow \mathcal{H}om(Q,E) \longrightarrow \mathcal{H}om(Q,F) \longrightarrow \mathcal{H}om(Q,Q) \longrightarrow
\end{equation}
$$
\longrightarrow \mathcal{E}xt^{1}(Q,E) \longrightarrow \mathcal{E}xt^{1}(Q,F) \longrightarrow \mathcal{E}xt^{1}(Q,Q) \longrightarrow
$$
$$
\longrightarrow \mathcal{E}xt^{2}(Q,E) \longrightarrow \mathcal{E}xt^{2}(Q,F) \overset{\phi}{\longrightarrow} \mathcal{E}xt^{2}(Q,Q) \longrightarrow \mathcal{E}xt^{3}(Q,E).
$$
Since the sheaves $E$ and $F$ are torsion-free we have that $\mathcal{H}om(Q,E) = \mathcal{H}om(Q,F) = 0$. Note that the smooth curve $C$ and $0$-dimensional subscheme $W$ are locally complete intersections, so for any point $x \in \mathbb{P}^{3}$ we have the following
\begin{equation}\label{local_isom}
\text{Ext}^1_{\mathcal{O}_{\mathbb{P}^{3},x}}(\mathcal{O}_{C,x},\mathcal{O}_{\mathbb{P}^{3},x}) = 0, \ \ \ \ \text{Ext}^{1, 2}_{\mathcal{O}_{\mathbb{P}^{3},x}}(\mathcal{O}_{W,x},\mathcal{O}_{\mathbb{P}^{3},x}) = 0.
\end{equation}
These equalities together with the condition (\ref{supports}) immediately imply that the sheaf $\mathcal{E}xt^{1}(Q,F)$ is equal to zero, so we have the isomorphism
\begin{equation}\label{isom1}
\mathcal{E}xt^{1}(Q,E) \simeq \mathcal{H}om(Q,Q).
\end{equation}
Also (\ref{local_isom}) and (\ref{supports}) imply that $\text{Supp}(\mathcal{E}xt^{2}(Q,F)) \subset C$, so we necessarily have the inclusion $\text{im}~\phi \subset \mathcal{E}xt^{2}(Q,Q)|_{C}$. On the other hand, homological dimension of the structure sheaf $\mathcal{O}_{C}$ is equal to 2, so we also have $\text{Supp}(\mathcal{E}xt^{3}(Q,E)) \cap C = \emptyset$. Now suppose that $\text{im}~\phi \subsetneq \mathcal{E}xt^{2}(Q,Q)|_{C}$, then $\text{Supp}(\text{coker}~\phi) \cap C \neq \emptyset$ which leads to the contradiction $\text{Supp}(\mathcal{E}xt^{3}(Q,E)) \cap C \neq \emptyset$ because $\text{coker}~\phi \hookrightarrow \mathcal{E}xt^{3}(Q,E)$. Also note that $\mathcal{E}xt^{2}(Q,Q) = \mathcal{E}xt^{2}(L,L) \oplus \mathcal{E}xt^{2}(\mathcal{O}_{W},\mathcal{O}_{W})$ due to $C \cap W = \emptyset$, so $\mathcal{E}xt^{2}(Q,Q)|_{C} = \mathcal{E}xt^{2}(L,L)$. Therefore, $\text{im}~\phi = \mathcal{E}xt^{2}(L,L)$ and we  have the following exact sequence
\begin{equation}\label{exact seq for ext2_0}
0 \longrightarrow \mathcal{E}xt^{1}(Q,Q) \longrightarrow \mathcal{E}xt^{2}(Q,E) \longrightarrow \mathcal{E}xt^{2}(L,F) \longrightarrow \mathcal{E}xt^{2}(L,L) \longrightarrow 0.
\end{equation}
Next, since $c_1(F)=0$ and $\text{Sing}(F) \cap C = \emptyset$, we have $\text{det}(F \otimes \mathcal{O}_{C}) \simeq \text{det}(F) \otimes \mathcal{O}_{C} \simeq \mathcal{O}_{\mathbb{P}^{3}} \otimes \mathcal{O}_{C} \simeq \mathcal{O}_{C}$. Therefore, the following exact triple holds
\begin{equation}\label{restriction of refl}
0 \longrightarrow L^{-1} \longrightarrow F \otimes \mathcal{O}_{C} \longrightarrow L \longrightarrow 0.
\end{equation}
Note that $\mathcal{E}xt^{2}(L, F) \simeq \mathcal{E}xt^{2}(L,F \otimes \mathcal{O}_{C})$ because $F$ is locally-free along $C$. In particular, it means that $\mathcal{E}xt^{2}(L, F)$ is locally-free $\mathcal{O}_{C}$-sheaf. Since $\mathcal{E}xt^{2}(L,L^{-1})$ is also locally-free $\mathcal{O}_{C}$-sheaf of rank 1, then applying the functor $\mathcal{H}om(L,-)$ to the triple (\ref{restriction of refl}) we obtain the following exact triple
\begin{equation}\label{ex_tr_1}
0 \longrightarrow \mathcal{E}xt^{2}(L,L^{-1}) \longrightarrow \mathcal{E}xt^{2}(L,F) \longrightarrow \mathcal{E}xt^{2}(L,L) \longrightarrow 0.
\end{equation}
Moreover, we have the following commutative diagram
\begin{equation*}
\begin{tikzcd}[column sep=small]
\mathcal{E}xt^{2}(L,F)  \rar \dar{\simeq} & \mathcal{E}xt^{2}(L,L) \dar{=} \\
\mathcal{E}xt^{2}(L,F \otimes \mathcal{O}_{C})  \rar & \mathcal{E}xt^{2}(L,L)
\end{tikzcd}
\end{equation*}
So the morphism $\mathcal{E}xt^{2}(L,F) \longrightarrow \mathcal{E}xt^{2}(L,L)$ in the triple (\ref{ex_tr_1}) coincides with the last morphism in the exact sequence (\ref{exact seq for ext2_0}). Therefore, we can simplify (\ref{exact seq for ext2_0}) as 
\begin{equation}\label{exact seq for ext2}
0 \longrightarrow \mathcal{E}xt^{1}(Q,Q) \longrightarrow \mathcal{E}xt^{2}(Q,E) \longrightarrow \mathcal{E}xt^{2}(L,L^{-1}) \longrightarrow 0.
\end{equation}
Note that for any subscheme $X \subset \mathbb{P}^{3}$ we have that 
$$
\mathcal{E}xt^{1}(\mathcal{O}_{X}, \mathcal{O}_{X}) \simeq \mathcal{H}om(I_{X}, \mathcal{O}_{X}) \simeq \mathcal{H}om(I_{X}/I_{X}^{2}, \mathcal{O}_{X}) = N_{X/\mathbb{P}^{3}}.
$$ 
(the last equality is the definition of the normal sheaf). Besides, if $X$ is a locally complete intersection of the pure dimension $1$ then (see \cite[Prop. 7.5]{AG})
$$
\mathcal{E}xt^{2}(\mathcal{O}_{X}, \mathcal{O}_{X}) \simeq \mathcal{E}xt^{2}(\mathcal{O}_{X}, \mathcal{O}_{\mathbb{P}^{3}}) \simeq \mathcal{E}xt^{2}(\mathcal{O}_{X}, \omega_{\mathbb{P}^{3}})(4) \simeq \omega_{X}(4).
$$
Now consider the case $X = C \sqcup W, \ Q = L \oplus \mathcal{O}_{W}$. Since $L$ is an invertible $\mathcal{O}_{C}$-sheaf it follows that 
$$
\mathcal{E}xt^{1}(L,L) \simeq \mathcal{E}xt^{1}(\mathcal{O}_{C},\mathcal{O}_{C}), \ \ \ \mathcal{E}xt^{2}(L,L^{-1}) \simeq \mathcal{E}xt^{2}(\mathcal{O}_{C}, \mathcal{O}_{C}) \otimes L^{-2}.
$$
From these formulas one can deduce the isomorphisms
$$
\mathcal{E}xt^{1}(Q,Q) \simeq N_{C / \mathbb{P}^{3}} \oplus N_{W / \mathbb{P}^{3}}, \ \ \ \mathcal{E}xt^{2}(L,L^{-1}) \simeq \omega_{C}(4) \otimes L^{-2}.
$$
Substituting them to (\ref{exact seq for ext2}) we obtain the following exact triple
\begin{equation}\label{triple1 for ext1}
0 \longrightarrow N_{C / \mathbb{P}^{3}} \oplus N_{W / \mathbb{P}^{3}} \longrightarrow \mathcal{E}xt^{2}(Q,E) \longrightarrow \omega_{C}(4) \otimes L^{-2} \longrightarrow 0.
\end{equation}

After applying the functor $\mathcal{H}om(-,F)$ to the exact triple (\ref{main}) we obtain the long exact sequence of sheaves
\begin{equation}\label{v2}
0 \longrightarrow \mathcal{H}om(Q, F) \longrightarrow \mathcal{H}om(E,F) \longrightarrow \mathcal{H}om(F,F) \longrightarrow
\end{equation}
$$
\longrightarrow \mathcal{E}xt^{1}(Q, F) \longrightarrow \mathcal{E}xt^{1}(E, F) \longrightarrow \mathcal{E}xt^{1}(F, F) \longrightarrow
$$
$$
\longrightarrow \mathcal{E}xt^{2}(Q, F) \longrightarrow \mathcal{E}xt^{2}(E, F) \longrightarrow \mathcal{E}xt^{2}(F, F).
$$
As it was already explained $\mathcal{H}om(Q, F)=\mathcal{E}xt^{1}(Q, F)=0$ and the morphism $\mathcal{E}xt^{1}(F, F) \longrightarrow \mathcal{E}xt^{2}(Q, F)$ is zero, so we have the following isomorphisms
\begin{equation}\label{isom2}
\mathcal{H}om(E,F) \simeq \mathcal{H}om(F,F), \ \ \ \mathcal{E}xt^{1}(E, F) \simeq \mathcal{E}xt^{1}(F, F).
\end{equation}

Consider the part of the commutative diagram with exact rows and columns obtained by applying the bifunctor $\mathcal{H}om(-,-)$ and its derivative $\mathcal{E}xt(-,-)$ to the exact triple (\ref{main}) which looks as follows
\begin{equation*}
\begin{tikzcd}[column sep=small]
  & 0  \dar & 0 \dar &  & & \\
0 \rar & \mathcal{H}om(F,E) \rar \dar & \mathcal{H}om(F,F) \rar \dar & \mathcal{H}om(F,Q) \rar & 0 & \\
0 \rar & \mathcal{H}om(E,E) \rar{\tau} \dar & \mathcal{H}om(E,F) \dar & & & \\
 & \mathcal{E}xt^{1}(Q, E) \dar  & 0  & & & \\
  & 0  &  &
\end{tikzcd}
\end{equation*}
Note that the uppermost horizontal and the leftmost vertical exact triples of this commutative diagram coincide with the exact triples (\ref{triple1}) and (\ref{triple2}), respectively. Due to the isomorphism (\ref{isom2}) the sheaf $\text{coker}~\tau$ fits into the exact triple
\begin{equation}\label{coker1}
0 \longrightarrow \mathcal{H}om(E, E) \longrightarrow \mathcal{H}om(F, F) \longrightarrow \text{coker}~\tau \longrightarrow 0.
\end{equation}
On the other hand, applying the Snake Lemma to the commutative diagram above and using the isomorphism (\ref{isom1}) we have the exact triple
\begin{equation}\label{coker2}
0 \longrightarrow \mathcal{H}om(Q, Q) \longrightarrow \mathcal{H}om(F,Q) \longrightarrow \text{coker}~\tau \longrightarrow 0.
\end{equation}
Since $\text{dim}~W=0$ we have that $h^1(\mathcal{H}om(F,\mathcal{O}_{W})) = 0$. Therefore, the condition (\ref{cond for h1}) implies $h^1(\mathcal{H}om(F,Q)) = h^1(\mathcal{H}om(F,L)) = 0$, so from the triple (\ref{coker2}) we obtain that
$$
h^{0}(\text{coker}~\tau) = h^{0}(\mathcal{H}om(F,Q))-h^{0}(\mathcal{H}om(Q, Q)) + h^{1}(\mathcal{H}om(Q, Q)),
$$
$$
h^{1}(\text{coker}~\tau)=h^{2}(\text{coker}~\tau)=0.
$$
Using these equalities and the fact that the sheaf $E$ is simple due to its stability, the triple (\ref{coker1}) implies the following formula
\begin{equation}\label{h1}
h^{1}(\mathcal{H}om(E, E))=1-h^{0}(\mathcal{H}om(F, F)) + h^{0}(\mathcal{H}om(F,Q))-
\end{equation}
$$
-h^{0}(\mathcal{H}om(Q, Q)) + h^{1}(\mathcal{H}om(Q, Q)) + h^{1}(\mathcal{H}om(F, F))=
$$
$$
=\text{dim~}\text{Hom}(F, Q) / \text{Aut}(Q) - \text{dim }\text{P}\text{Aut}(F) + \text{dim Jac}(C) + h^{1}(\mathcal{H}om(F, F)),
$$
\begin{equation}\label{h2}
h^{2}(\mathcal{H}om(E,E))=h^{2}(\mathcal{H}om(F,F)).
\end{equation}
Next, using the isomorphisms (\ref{isom1 for ext1}), (\ref{isom2 for ext1}), the triple (\ref{triple1 for ext1}) and the condition (\ref{cond defect}) we obtain the formula
\begin{equation}\label{h0}
h^0(\mathcal{E}xt^{1}(E,E)) = h^{0}(N_{W/\mathbb{P}^{3}}) + h^{0}(N_{C/\mathbb{P}^{3}}) + h^0(\mathcal{E}xt^{1}(F,F))=
\end{equation}
$$
=\text{dim}~\mathcal{H}_{0} \times \mathcal{H}_{1} + h^0(\mathcal{E}xt^{1}(F,F)).
$$
Substituting the formulas (\ref{h1}-\ref{h0}) to the equality (\ref{equality for ext1}) and using (\ref{equality for ext1f}), we obtain the following formula
\begin{equation}
\text{dim Ext}^1(E,E)=\text{dim~}\mathcal{R} + \text{dim} \ \mathcal{H}_{0} + \text{dim} \ \mathcal{H}_{1} + \text{dim Jac}(C) +
\end{equation}
$$
+ \text{dim~}\text{Hom}(F, Q) / \text{Aut}(Q) - \text{dim }\text{P}\text{Aut}(F).
$$
Now taking into account (\ref{picard_dim}) and (\ref{dim}) we immediately obtain the statement of the theorem.
\hfill$\Box$

\vspace{3mm}

From the construction of the component $\overline{\mathcal{C}}$ it follows that the general sheaf $E$ of $\overline{\mathcal{C}}$ has singularities of mixed dimension, more precisely, we have $\text{Sing}(E)=C \sqcup \text{Sing}(E^{\vee \vee}) \sqcup W$, where $C$ is a curve of degree more than $1$, $\text{dim} \ \text{Sing}(E^{\vee \vee}) = \text{dim} \ W = 0$ and $\text{Sing}(E^{\vee \vee}) \neq \emptyset$. On the other hand, the general sheaves of all known components parameterizing sheaves with mixed singularities have singularity sets of the form $l \sqcup W$, where $l$ is a projective line and $\text{dim}~W = 0$. Therefore, the component $\overline{\mathcal{C}}$ is not one of the previously known components.

Next, note that the construction of the open subset $\mathcal{C}$ and its closure $\overline{\mathcal{C}} \subset \mathcal{M}(m+d)$ depends on the choice of the number $s$ of disjoint points, the choice of the component $\mathcal{R}$ from two series $\mathcal{S}_{a,b,c}, \ \mathcal{V}_{m}$, and the choice of the Hilbert scheme $\mathcal{H}_{1}$ from two series of the Hilbert schemes $\text{Hilb}_{d}, \ \text{Hilb}_{(d_1, d_2)}$. So, in fact, we have the series of components which we will denote by $\overline{\mathcal{C}(\mathcal{R}, \mathcal{H}_{1}, s)}$. 

Also it is worth to note that the described series of components can be extended to a larger series of components. In order to construct them we consider strictly $\mu$-semistable reflexive sheaves defined by the triple (\ref{serre for non-stable}), where $Y$ is a disjoint union of rational curves. Then we do elementary transformations of these reflexive sheaves along disjoint union of a collection of distinct points, smooth rational curves and complete intersection curves, simultaneously. It seems that the proof of irreducibility of components of this extended series is essentially the same as above and need only minor modifications.

Since a complete enumeration of components of $\mathcal{M}(k)$ for small values of $k$ is of particular interest, we point out that the series described above contains a new component from $\mathcal{M}(3)$, namely, $\overline{\mathcal{C}(\mathcal{V}_{1}, \text{Gr}(2,4),0)}$. By construction the general sheaf $[E]$ of this component fits into the exact sequence
$$
0 \longrightarrow E \longrightarrow F \longrightarrow \mathcal{O}_{C}(2) \longrightarrow 0,
$$
where $C$ is a smooth conic and $[F] \in \mathcal{V}_{1} = \mathcal{V}(0,1,2)$ is a $\mu$-semistable sheaf satisfying the following exact triple
$$
0 \longrightarrow \mathcal{O}_{\mathbb{P}^{3}} \longrightarrow F \longrightarrow I_{l} \longrightarrow 0,  \ \ \ l \in \text{Gr}(2,4).
$$
Dimension of the component $\mathcal{C}(\mathcal{V}_{1}, \text{Gr}(2,4),0)$ is equal to $21$ and its spectrum is $(-1,0,1)$. Therefore, the number of components of $\mathcal{M}(3)$ is at least 11.

\end{document}